\documentclass[a4paper,11pt,english,spanish,leqno]{amsart}
\usepackage[utf8x]{inputenc}
\usepackage[bitstream-charter]{mathdesign}
\usepackage{amsmath, xcolor,amsthm,epsfig,epstopdf,url,array}
\usepackage[colorlinks=true]{hyperref}
\usepackage[all]{xy}
\usepackage{tikz}
\tikzset{node distance=2cm, auto}
\usepackage{subcaption}
\usepackage[backrefs]{amsrefs}

\theoremstyle{plain}
\newtheorem{thm}{Theorem}[section]

\newtheorem{lem}[thm]{Lemma}

\newtheorem{prop}[thm]{Proposition}

\newtheorem{cor}[thm]{Corollary}

\newtheorem{obs}[thm]{Remark}
\newtheorem{problem}[thm]{Problem}

\theoremstyle{definition}
\newtheorem{defn}[thm]{Definition}

\theoremstyle{remark}
\newtheorem*{dem}{Proof}

\theoremstyle{plain}

\def\S#1{\mathbb{S}^{#1}}
\def\Esp#1#2{\text{\rm E}_{#1}\left[#2\right]}

\def\E#1{\mathcal{E}_{#1}}

\def\Xr{\mathfrak{X}_{*}^{(r)}}
\def\ra{\rightarrow}
\def\R{\mathbb{R}}
\def\C{\mathbb{C}}
\def\i{\mathbf{i}}
\def\M{\mathrm{M}}
\def\P{\mathbb{P}(\mathbb{C}^{d+1})}

\begin{document}

\title{The projective ensemble and distribution of points in odd--dimensional spheres}
\author{Carlos Beltr\'an and Uju\'e Etayo}
\date{\today{}}
\thanks{This research has been partially supported by Ministerio de Economía y Competitividad, Gobierno de España, through grants MTM2014-57590-P and MTM2015-68805-REDT.}

\begin{abstract}

We define a determinantal point process on the complex projective space that reduces to the so--called spherical ensemble for complex dimension $1$ under identification of the $2$--sphere with the Riemann sphere. 
Through this determinantal point process we propose a point processs in odd-dimensional spheres that produces fairly well--distributed points, in the sense that the expected value of the Riesz $2$--energy for these collections of points is smaller than all previously known bounds.

\end{abstract}

\maketitle



\section{Introduction}

Given $s\in(0,\infty)$, the Riesz $s$--energy of a set on points $\omega_{n} = \{ x_{1}, \ldots  ,x_{n} \}$ on a subset $X\subseteq\R^{m}$ is
\begin{equation}\label{eq:Riesz}
\E{s}(\omega_{n}) = \displaystyle\sum_{i \neq j} \frac{1}{\| x_{i} - x_{j} \|^{s}}.
\end{equation} 
This energy has a physical interpretation for some particular values of $s$, i.e. for $s=1$ the Riesz energy is the Coulomb potential and for $s=d-2$ ($d\geq 3$) is the Newtonian potential. In the special case $s=0$ the energy is defined by 
\[
\E{0}(\omega_{n})=\left. \frac{d}{ds} \right|_{s=0}\E{s}(\omega_n) = \sum_{i \neq j} \log\| x_{i} - x_{j} \|^{-1} 
\]
and is related to the transfinite diameter and the capacity of the set by classical potential theory, see for example \cite{doohovskoy2011foundations}.

The minimal value of this energy and its asymptotic behavior have been extensively studied, most remarkably in the case that $X=\mathbb{S}^{d}\subseteq\R^{d+1}$ is the $d$--dimensional unit sphere.
In \cite{10.2307/117605} it was proved that for $d > 2$ and $0<s<d$ there exist constants $c>C> 0$ (depending only on $d$ and $s$) such that
\begin{equation}\label{eq:bounds}
-cn^{1 + \frac{s}{d}} 
\leq
\min_{\omega_n}\left(\E{s}(\omega_{n})\right) - V_{s}(\mathbb{S}^{d})n^{2} 
\leq 
-Cn^{1 + \frac{s}{d}},
\end{equation}
where $V_{s}(\mathbb{S}^{d})$ is the continuous s-energy for the normalized Lebesgue measure,
\begin{equation}
V_{s}(\mathbb{S}^{d})
=
\frac{1}{Vol(\mathbb{S}^{d})^{2}} \displaystyle\int_{p,q \in \mathbb{S}^{d}} \frac{1}{\|p-q\|^{s}}dpdq
=
2^{d-s-1} \frac{\Gamma \left( \frac{d+1}{2} \right) \Gamma \left( \frac{d-s}{2} \right)}{\sqrt{\pi} \Gamma \left(  d- \frac{s}{2} \right)}.
\end{equation}
Finding the precise value of the constants in \eqref{eq:bounds} is an important open problem and has been addressed by several authors, see \cites{ BHS2012b, Sandi, LB2015, MR1306011} for some very precise conjectures and \cite{Brauchart2015293} or \cite{BHSlibro} for surveys. One can post the problem as follows
\begin{problem}
 For $s\in(0,d)$, let $C_{s,d,n}$ be defined by
 \[
  0<C_{s,d,n}=\frac{ V_{s}(\mathbb{S}^{d})n^{2} -\min_{\omega_n}\left(\E{s}(\omega_{n})\right)}{n^{1+\frac{s}{d}}}.
 \]
 Find asymptotic values for $C_{s,d,n}$ as $n\to\infty$. In particular, prove if the limit exists.

\end{problem}
A sometimes successful strategy for the upper bound in the constant $C_{s,d,n}$ is to take collections of random points in $\S{d}$ and then compute the expected value of the energy (which is of course greater than or equal to the minimum possible value). Simply taking $n$ points with the uniform distribution in $\S{d}$ already gives the correct term $V_{s}(\mathbb{S}^{d})n^{2}$, and other distributions with nice separation properties have proved successful in bounding the constant $C_{s,d,n}$. 

We are thus interested in computationally feasible random procedures to generate points in sets which exhibit local repulsion. One natural choice is using {\em determinantal point processes} which have these two properties (see \cite{Hough_zerosof} for theoretical properties and \cite{PhysRevE.79.041108} for an implementation). A brief summary of the fundamental properties of determinantal point processes is given in Section \ref{sec:DPP}.

In a recent paper \cite{EJP3733} a determinantal point process named the {\em spherical ensemble} is used to produce low--energy random configurations on $\S{2}$. This process was previously studied by Krishnapur \cite{krishnapur2009} who proved a remarkable fact: the spherical ensemble is equivalent to taking eigenvalues of $A^{-1}B$ (where $A,B$ have Gaussian entries) and sending them to the sphere through the stereographic projection.

In \cite{BMOC2015energy} a different determinantal point process rooted on the use of spherical harmonics is described, producing low--energy random configurations in $\S{d}$ for some infinite sequence of values of $n$. In particular, it is proved in that paper that
\begin{equation}\label{eq:BMOCliminf}
 \limsup_{n\to\infty} C_{s,d,n}\geq 2^{s-s/d}V_s(\S{d})
\frac{ d\,\Gamma\left( 1+\frac{d}{2} \right)   \Gamma\left( \frac{1+s}{2} 
\right)\Gamma\left(d-\frac{s}{2}\right) }
{ \sqrt{\pi} \Gamma\left( 1+\frac{s}{2} \right) \Gamma\left( 1+\frac{s+d}{2} 
\right)\left(d!\right)^{1-\frac{s}{d}}},\quad 0<s<d.
\end{equation}
If $d-1<s<d$, $\limsup$ can be changed to $\liminf$ in \eqref{eq:BMOCliminf} (see \cite[Cor. 2]{BMOC2015energy}).
 The bound in \eqref{eq:BMOCliminf} is the best known to the date for general $d$ (although more precise bounds exist for particular values of $d$ including $d=1,2$, see \cites{BLMS:BLMS0621, LB2015}). In particular, for  $s=2$ and odd dimensions the formula in \eqref{eq:BMOCliminf} reads
\begin{equation}\label{eq:BMOCliminf2}
 \limsup_{n\to\infty} C_{2,2d+1,n}\geq \frac{2^{1-\frac{2}{2d+1}}\left({(2d+1)!}\right)^\frac{2}{2d+1}}{(2d-1)(2d+3)}\stackrel{d\to\infty}{\to}\frac{2}{e^2}.
\end{equation}

The determinantal point process in \cite{BMOC2015energy} is called the {\em harmonic ensemble} and it is shown to be the optimal one (at least for $s=2$) among a certain class of determinantal point processes obtained from subspaces of functions with real values defined in $\S{d}$.

However, the bound in \cite{BMOC2015energy} for the case $d=2$ is worse than that of \cite{EJP3733}, which is not surprising since the spherical ensemble uses complex functions and is thus of a different nature. 

 An alternative natural interpretation of Krishnapur's result is to consider eigenvalues $(\alpha,\beta)\in\mathbb P(\C^2)$ of the generalized eigenvalue problem $\det(\beta A-\alpha B)=0$ and to identify $\mathbb P(\C^2)$ with the Riemann sphere. An homotety then generates the points in the unit sphere $\S{2}$. This remark suggests that the spherical ensemble can be seen as a natural point process in the complex projective space, and a search for an extension to higher dimensions is in order. In this paper we extend this process in a very natural manner to $\P$ for any $d\geq2$. We will propose the name {\em projective ensemble}.

In order to show the separation properties of the projective ensemble we will define a (probably non--determinantal) point process in odd--dimensional spheres, which will allow us to compare our results to those of \cite{BMOC2015energy}. This point process is as follows: first, choose a number $r$ of points in $\P$ coming from the projective ensemble. Then, consider $k$ equally spaced unit norm affine representatives of each of the projective points. We allow these points to be rotated by a randomly chosen phase. As a result, we get $rk$ points in the odd--dimensional sphere $\mathbb{S}(\C^{d+1})\equiv\mathbb{S}^{2d+1}$. 
 
We study the expected $2$--energy of such a point process. Our first main result can be succinctly written as follows.

\begin{figure}[htp]
	\centering
		\includegraphics[width=0.9\textwidth]{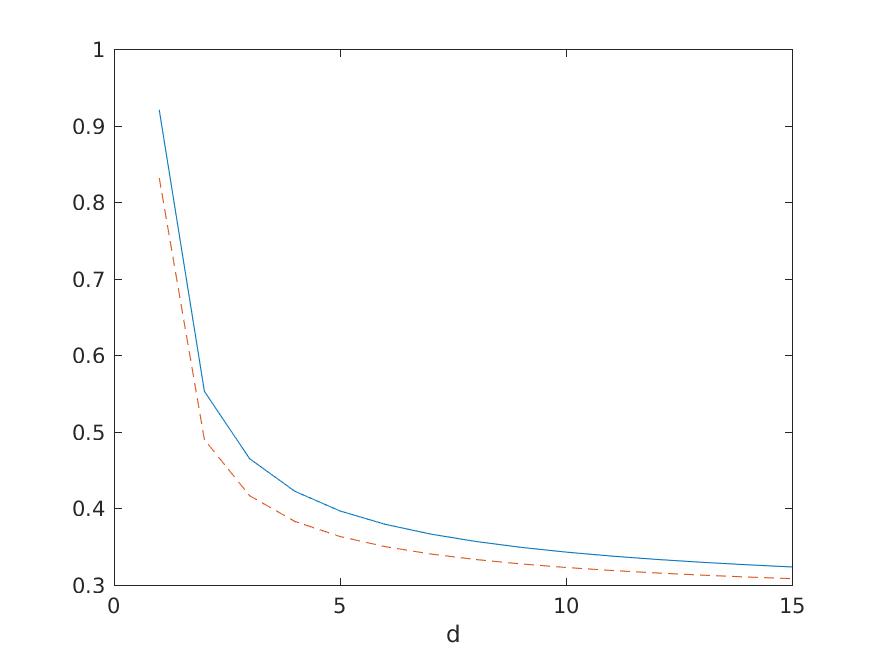}
	\caption{Comparison of the values of the constants in \eqref{eq:nuestracota} (blue solid line) and \eqref{eq:BMOCliminf2} (red dashed line).} \label{fig:comparacion}
\end{figure}

\begin{thm}\label{th:main}
With the notations above,

\begin{equation}\label{eq:nuestracota}
\begin{split}
&\limsup_{n\to\infty}  C_{2,2d+1,n}
\geq 
\frac{3^{1 - \frac{2}{2d+1}} (2d-1)^{1 - \frac{2}{2d+1}} (2d+1) \Gamma\left( d- \frac{1}{2} \right)^{2 - \frac{2}{2d+1}}}{2^{4 - \frac{2}{2d+1}} (d!)^{2 - \frac{4}{2d+1}}}\stackrel{d\to\infty}{\to}\frac{3}{4e}.
\end{split}
\end{equation}
\end{thm}
The bound in Theorem \ref{th:main} is larger than that of \eqref{eq:BMOCliminf2}, which shows that random configurations of points coming from this point process are, at least from the point of view of the $2$--energy, better distributed than those coming from the harmonic ensemble. See Figure \ref{fig:comparacion} for a graphical comparison of both bounds.

Since the point process we have defined in $\mathbb{S}^{2d+1}$ starts by choosing points in $\P$ coming from the projective ensemble, Theorem \ref{th:main} gives us arguments to think that the projective ensemble produces quite well distributed points in $\P$ (for $d=1$ this property is quantitatively described in \cite{EJP3733}). There are several ways to measure how well distributed a collection of points is in $\P$. For example, one can study the natural analogues of Riesz's energy as in Theorem \ref{Prop_proy} below. A very natural measure is given by the value of {\em Green's energy} of \cite{Juan}: let $G:\P\times\P\to[0,\infty]$ be the Green function of $\P$, that is, $G(x,\cdot)$ is zero--mean for all $x$, $G$ is symmetric and $\Delta_yG(x,y)=\delta_x(y)-vol(\P)^{-1}$, with $\delta_x$ the Dirac's delta function, in the distributional sense. The Green energy of a collection of $r$ points $\omega_r=(x_1,\ldots,x_r)\in\P^r$ is defined as
 \[
  \E{G}(\omega_r)=\sum_{i\neq j}G(x_i,x_j).
  \]
    Minimizers of Green's energy are assymptotically well--distributed (see \cite[Main Theorem]{Juan}). Our second main result will follow from the computation of the expected value of Green's energy for the projective ensemble.
\begin{thm}\label{th:main2}
 Let $d\geq 2$. Then,
 \begin{equation}\label{eq:final}
  \liminf_{r\to\infty}\frac{\min_{\omega_r}\left(\E{G}^{\mathbb{P}}(\omega_{r})\right)}{r^{2 - \frac{1}{d}}}\leq-\frac{ (d!)^{1 - \frac{1}{d}}}{4\pi^{d} (d-1)}.
 \end{equation}

\end{thm}
Theorem \ref{th:main2} gives a criterium to decide how well--distributed a collection of projective points is: compute their Green's energy and compare to that of \eqref{eq:final}.
%

%
%
%

\section{Determinantal point processes}\label{sec:DPP}
\subsection{Basic notions}
In this section we follow \cite{Hough_zerosof}.

\begin{defn}\label{def:pp}
Let $\Lambda$ be a locally compact, polish topological space with a Radon measure $\mu$. A \textit{simple point process} $\mathfrak{X}$ of $n$ points in $\Lambda$ is a random variable taking values in the space of $n$ point subsets of $\Lambda$.
\end{defn}
There are some subtle issues in the general definition of point processes, see \cite[Section 1.2]{Hough_zerosof}. For our purposes we will only use simple point processes with a fixed, finite number of points.

\noindent For some point processes there exist \textit{joint intensities} satisfying the following definition.
\begin{defn}
 Let $\Lambda,\mathfrak{X}$ be as in Definition \ref{def:pp}. The joint intensities are functions (if any exist) $\rho_{k}:\Lambda^k\ra[0,\infty)$, $k\geq1$ such that for any family of mutually disjoint subsets $D_1,\ldots,D_k$ of $\Lambda$ we have
 \[
  \Esp{x\sim\mathfrak{X}}{\left(\prod_{i=1}^k \sharp(x\cap D_i)\right)}=\int_{\prod D_i}\rho_k(x_1,\ldots,x_k)\,d\mu(x_1,\ldots,x_k).
 \]
Here, $\mathrm{E}$ denotes expectation and by $x\sim\mathfrak{X}$ we mean that $x$ is a subset of $\Lambda$ with $n$ elements, obtained from the point process $\mathfrak{X}$.
\end{defn}
From \cite[Formula (1.2.2)]{Hough_zerosof}, for any measurable function $\phi: \Lambda^{k} \longrightarrow [0, \infty)$ the following equality holds.
	\begin{multline}
	\label{eq_1}
	\Esp{x\sim\mathfrak{X}}{\sum_{i_{1} \ldots  i_{k}\text{ distinct}} \phi(x_{i_{1}} ,\ldots , x_{i_{k}})} \\ =\int_{y_{1},\ldots ,y_{k} \in \Lambda} \phi(y_{1},\ldots ,y_{k}) \rho_{k}(y_{1},\ldots ,y_{k})\,d\mu(y_1,\ldots,y_k).
	\end{multline}
	
\noindent  Sometimes these intensity joint functions can be written as $\rho_{k}(x_{1},\ldots ,x_{k}) = det(K(x_{i}, x_{j})_{i,j = 1,\ldots ,k})$ for some function $K:\Lambda\times\Lambda\to\C$.
In this case, we say that $\mathfrak{X}$ is a  \textit{determinantal point process}. 
A particularly amenable collection of such processes is obtained from $n$--dimensional subspaces of the Hilbert space $\mathrm{L}^2(\Lambda,\C)$ (i.e. the set of square--integrable complex functions in $\Lambda$). 
Recall that the reproducing kernel of $H$ is the unique continuous, skew--symmetric, positive--definite function $K_H:\Lambda\times\Lambda\to\C$ such that
\[
f(x)=\langle f,K_H(\cdot,x)\rangle=\int_{y\in\Lambda}f(y)K_H(x,y)\,dy,\quad x\in\Lambda,f\in H.
\]
Given any orthonormal basis $\varphi_{1},\ldots ,\varphi_{n}$ of $H$, we have
\begin{equation}\label{eq:kernel}
 K_H(x,y)=\sum_{i=1}^n\varphi_i(x)\overline{\varphi_i(y)}.
\end{equation}
Such a kernel $K_H$ is usually called a {\em projection kernel} of trace $n$.	
\begin{prop}\label{prop:MS}	
Let $\Lambda$ be as in Definition \ref{def:pp} and let $H \subset \text{\rm L}^{2}(\Lambda,\C)$ have dimension $n$. Then there exists a point process $\mathfrak{X}_H$ in $\Lambda$ of $n$ points with associated join intensity functions
	\begin{equation*}
	\rho_{k}(x_{1},\ldots ,x_{k}) = det(K_H(x_{i}, x_{j})_{i,j = 1,\ldots ,k}).
	\end{equation*}
In particular for any measurable function $f:\Lambda \times \Lambda \longrightarrow [0,\infty)$ we have
	\begin{multline}
	\label{th_1}
	\Esp{x \sim \mathfrak{X}_H}{
	\displaystyle\sum_{i \neq j} f(x_{i}, x_{j})}  
	\\
	= \int_{p,q \in \Lambda}
	\left( K_H(p,p)K_H(q,q) - |K_H(p,q)|^{2} \right) f(p,q)\,d\mu(p,q).
	\end{multline}
We will call $\mathfrak{X}_H$ a projection determinantal point process with kernel $K_H$.
\end{prop}

\begin{dem}
This proposition is a direct consequence of the Macchi--Soshnikov Theorem, see \cite{Macchi, Soshni} or  \cite[Theorem 4.5.5]{Hough_zerosof}. $\square$
\end{dem}

\begin{obs}\label{obs:constant}
We note that in the hypotheses of Proposition \ref{prop:MS}, from \eqref{eq_1} with $\phi\equiv1$ we also have
\begin{equation*}
n=\Esp{x\sim\mathfrak{X}_H}{n}  
= \int_{p\in\Lambda} K_H(p,p)\,d\mu(p),
\end{equation*}

\noindent In particular, if $K_H(p,p)$ is constant then we must have $K_H(p,p)=n/Vol(\Lambda)$.
\end{obs}
\subsection{Transformation under diffeomorphisms}

We now describe the push--forward of a projection determinantal point process. We are most interested in the case that the spaces are Riemannian manifolds (which are locally compact, Polish and measurable spaces).

\begin{prop}
\label{prop_1}
Let $\M_{1}$ and $\M_{2}$ be two Riemannian manifolds and let $\phi: M_{1} \longrightarrow \M_{2}$ be a $\boldsymbol{\mathcal{C}^{1}}$ diffeomorphism. Let $H \subset \mathrm{L}^{2}(\M_{1},\C)$ be an $n$--dimensional subspace. Then, the set
\begin{align*}
	H_{*} =& \left \{ f: \M_{2} \longrightarrow \mathbb{C}: \sqrt{|\text{\rm Jac}(\phi)(x)|}( f \circ \phi)(x) \in H \right \} \\
	 =&\left \{ g \circ \phi^{-1}(\cdot)\sqrt{|\text{\rm Jac}(\phi^{-1})(\cdot)|}: g \in H \right \}
\end{align*}
is an $n$--dimensional subspace of $\mathrm{L}^2(\M_2,\C)$. Its associated determinantal point process $\mathfrak{X}_{H_*}$ has kernel
\begin{equation}
	\begin{split}
	K_{H_*}(a,b) &= \frac{K_H(\phi^{-1}(a), \phi^{-1}(b))}{\sqrt{|\text{\rm Jac}(\phi)(\phi^{-1}(a)) \text{\rm Jac}(\phi)(\phi^{-1}(b))|}}  \\
& = K_H(\phi^{-1}(a), \phi^{-1}(b))\sqrt{|\text{\rm Jac}(\phi^{-1})(a) \text{\rm Jac}(\phi^{-1})(b)|}.
	\end{split}
	\end{equation}
	(We are denoting by $\mathrm{Jac}$ the Jacobian determinant).
\end{prop}
This proposition is a direct consequence of the change of variables formula, see Section \ref{subsec:proof1} for a short proof.


\section{The projective ensemble}\label{sec:iran}

Consider the standard Fubini--Study metric in the complex projective space of complex dimension $d$, denoted by $\P$. The distance between two points $p,q\in\P$ is given by:
	\begin{equation*}
	\sin d_{\P}(p,q) 
	=
	\sqrt{1 - \frac{| \left\langle p,q \right\rangle |^{2}}{|| q ||^{2} || p ||^{2}}}
	=
	\sqrt{1 - \left| \left\langle \frac{p}{||p||},\frac{q}{||q||}\right\rangle \right|^{2}}.
	\end{equation*}

\begin{defn}
\label{def_2}
\noindent Let $L \geq 0$ and consider the set of the following functions defined in $\mathbb{C}^{d}$:
\begin{equation}
\label{eq_3}
I_{d,L} = \left\lbrace \sqrt{C_{\alpha_{1},\ldots ,\alpha_{d}}^{L}} \frac{z_{1}^{\alpha_{1}}\ldots z_{d}^{\alpha_{d}}}{(1+\|z\|^{2})^{\frac{d+L+1}{2}}} \right\rbrace_{\alpha_{1}+\ldots +\alpha_{d} \leq L}
\end{equation}

\noindent where $\alpha_{1},\ldots ,\alpha_{d}$ are non--negative integers and

\begin{equation*}
C_{\alpha_{1},\ldots ,\alpha_{d}}^{L} 
= 
\frac{1}{{\pi}^{d}} \frac{(d+L)!}{\alpha_{1}!\ldots \alpha_{d}! (L-(\alpha_{1}+\ldots +\alpha_{d}))!}.
\end{equation*}

\noindent Let $\mathcal{H}_{d,L}=Span(I_{d,L})\subset\text{\rm L}^{2}(\mathbb{C}^{d},\C)$ which is a subspace of complex dimension $r = {d+L \choose d}$. 
The collection $I_{d,L}$ given in \eqref{eq_3} is an orthonormal basis of $\mathcal{H}_{d,L}$ (for the usual Lebesgue measure in $\mathbb{C}^{d}$) and the reproducing kernel $K: \mathbb{C}^{d} \times  \mathbb{C}^{d} \longrightarrow \mathbb{C}$ is given by:
\begin{equation*}
K(z,w) = \frac{r d!}{{\pi}^{d}} \frac{(1 + \left\langle z,w \right\rangle)^{L}}{(1+\|z\|^{2})^{\frac{d+L+1}{2}} (1+\|w\|^{2})^{\frac{d+L+1}{2}}}.
\end{equation*}
From Proposition \ref{prop:MS}, there is an associated determinantal point process of $r$ points in $\C^d$ that we denote by $\mathfrak{X}^{(r,d)}$.
\end{defn}

\begin{lem}\label{propdef_1}
Let $d\geq1$ and let $r$ be of the form $r={d+L \choose d}$. Then, the pushforward $\mathfrak{X}_*^{(r,d)}$ of $\mathfrak{X}^{(r,d)}$ under the mapping
\[
 \begin{matrix}
  \psi_d:&\C^d&\to&\P\\
  &z&\mapsto&(1,z)
 \end{matrix}
\]
is a determinantal point process in $\P$ whose associated kernel satisfies
\begin{equation*}
 \left| K_*^{(r,d)}(p,q) \right|
=
\frac{r d!}{\pi^{d}} 
\left| \left\langle \frac{p}{||p||}, \frac{q}{||q||} \right\rangle \right|^L.
\end{equation*}
We call this process the projective ensemble.
\end{lem}

\noindent See Section \ref{propLemma32} for a proof of Lemma \ref{propdef_1}.
The spherical ensemble described in \cite{krishnapur2009,EJP3733} is just the case $d=1$ of the projective ensemble identifying $\mathbb{P}(\mathbb{C}^{2})$ with the Riemmann sphere and translating the process to the unit sphere.

The next result computes the expected value of a Riesz--like energy for the projective ensemble.

\begin{thm}
\label{Prop_proy}
Let $L\geq1$. For $r=\binom{d+L}{d}$ and $\omega_{r}=(x_1,\ldots,x_r)\in\P^r$ let
\begin{equation*}
\mathcal{E}_{s}^{\mathbb{P}}(\omega_{r})
=
\sum_{i \neq j} \frac{1}{\sin \left( d_{\P} (x_{i}, x_{j}) \right)^{s}}.
\end{equation*}
Then, for $0 < s < 2d$,
\begin{equation*}
\begin{split}
 \Esp{\mathfrak{X}_{*}^{(r)}}{\mathcal{E}_{s}^{\mathbb{P}}(\omega_{r}) } 
& =
\frac{d}{d - \frac{s}{2}}r^{2} - r^{2}dB \left( d - \frac{s}{2}, L + 1 \right)\\&= \frac{d}{d - \frac{s}{2}}r^{2}-r^{1+\frac{s}{2d}}\frac{d \Gamma\left(d - \frac{s}{2}\right)}{(d!)^{1-\frac{s}{2d}}}+o\left(r^{1+\frac{s}{2d}}\right).
\end{split}
\end{equation*}
Note that $d/(d-s/2)$ is precisely the continuous $s$--energy for the uniform measure in $\P$.
\end{thm}
\begin{cor}\label{cor:logaritmica}
 Let $L\geq1$. For $r=\binom{d+L}{d}$ and $\omega_{r}=(x_1,\ldots,x_r)\in\P^r$ let
\begin{equation*}
\mathcal{E}_{0}^{\mathbb{P}}(\omega_{r})
=
\sum_{i \neq j} \log\frac{1}{\sin \left( d_{\P} (x_{i}, x_{j}) \right)}.
\end{equation*}
Then,
\begin{equation*}
\begin{split}
 \Esp{\mathfrak{X}_{*}^{(r)}}{\mathcal{E}_{0}^{\mathbb{P}}(\omega_{r}) } 
& =
\frac{r^2}{2d}+\frac{r^2d}{2}B(d,L+1)\sum_{j=0}^L\frac{1}{d+j}\\&= \frac{r^2}{2d}+\frac{r\log r}{2d}+o\left(r\log r\right).
\end{split}
\end{equation*}
\end{cor}
Theorem \ref{Prop_proy} and Corollary  \ref{cor:logaritmica} are proved in sections \ref{sec:proofofProp_proy} and \ref{sec:prooflog}.

\section{A new point process in odd-dimensional spheres}\label{section4}
 We now describe a point process of $n$ points, for certain values of $n$, in $\S{2d+1}$ in the following manner.
\begin{defn}\label{def:pps2dmas1}
 Given integers $d,\,k,\,L\geq0$, let $r=\binom{d+L}{d}$ and $n=k\,r$. We define the following point process of $n$ points in $\S{2d+1}$. First, let
 \[
  x_1,\ldots,x_r\in\P
 \]
 be chosen from the projective ensemble $\mathfrak{X}_*^{(r,d)}$. Choose, for each $x_i$, one affine representative (which we denote by the same letter). Then, let $\theta_1,\ldots,\theta_r\in[0,2\pi)$ be chosen uniformly and independently and define
 \begin{equation}\label{eq:yes}
  y_{i}^{j}
=
e^{\i\left( \theta_{i} + \frac{2\pi j}{k} \right)} x_{i},\quad 1\leq i\leq r,\;0\leq j\leq k-1. 
 \end{equation}
 We denote this point process by $\mathfrak{X}_{\S{2d+1}}^{(k,L)}$.
\end{defn}
Note that the way to generate a collection of $n$ points coming from $\mathfrak{X}_{\S{2d+1}}^{(k,L)}$ amounts to taking $r$ points from the projective ensemble and taking, for each of these points, $k$ affine unit norm representatives, uniformly spaced in the great circle corresponding to each point, with a random phase.

The following statement shows that the expected $2$--energy of points generated from the point process of Definition \ref{def:pps2dmas1} can be computed with high precision. It will be proved in Section \ref{sec:proofUju1}.

\begin{prop}\label{Uju1}
\begin{multline}
\label{eq3bis}
\Esp{\mathfrak{X}_{\S{2d+1}}^{(k,L)}}{\E{2}(y_{1}^{0},...,y_{1}^{k-1},y_{2}^{0},...,y_{r}^{k-1})}\\  =\frac{d}{2d-1}  (kr)^{2}+\frac{rk^{3}}{12}-\frac{d\Gamma\left(d-\frac12\right)}{2(d!)^{1-\frac1{2d}}} k^2r^{1+\frac{1}{2d}}+o\left( k^2r^{1+\frac{1}{2d}}\right).
\end{multline}
\end{prop}

Following the same ideas one can also compute the expected $s$-energy for $n$ points coming from the point process $\mathfrak{X}_{\S{2d+1}}^{(k,L)}$ for other even integer values $s \in 2\mathbb{N}$, and a bound can be found for other values of $s>0$. The computations, though, are quite involved.

Proposition \ref{Uju1} describes how different choices of $L$ (i.e. of $r$) and $k$ produce different values of the expected $2$--energy of the associated $n=rk$ points. An optimization argument is in order: for given $n$, which is the optimal choice of $r$ and $k$? Since we know from \eqref{eq:bounds} that the second order term in the assymptotics is $\sim n^{1+2/(2d+1)}=(rk)^{1+2/(2d+1)}$, it is easy to conclude that the optimal values of $r$ and $k$ satisfy:
\[
 k\sim r^{\frac{1}{2d}}.
\]
The following corollary follows inmediately from Proposition \ref{Uju1}. 
\begin{cor}\label{cor:casi}
If we choose $k=A r^{\frac{1}{2d}}$ for some $A\in\R$ making that quantity a positive integer, then:
\begin{multline}
\label{eq3bis2}
\Esp{\mathfrak{X}_{\S{2d+1}}^{(k,L)}}{\E{2}(y_{1}^{0},...,y_{1}^{k-1},y_{2}^{0},...,y_{r}^{k-1})}\\ = \frac{d}{2d-1} n^{2}
+
\left(
\frac{A^{2 - \frac{2}{2d+1}} }{12} 
-
\frac{d \Gamma\left( d - \frac{1}{2} \right) A^{1 - \frac{2}{2d+1}} }{2 \left(d!\right)^{1 - \frac{1}{2d}} }
\right) 
n^{1+\frac{2}{2d+1}}
+
o\left(n^{1+\frac{2}{2d+1}}\right).
\end{multline}
\end{cor}
The proof of our first main theorem will follow easily from Corollary \ref{cor:casi}.


\section{Proof of the main results}

\subsection{Proof of Proposition \ref{prop_1}}\label{subsec:proof1}

We first prove that $H_*\subseteq\mathrm{L}^2(\M_2,\C)$. Indeed, for $f\in H_*$ we have
\[
 \int_{y\in \M_2}|f|^2\,dy=\int_{y\in \M_2}|g\circ\phi^{-1}(y)|^2|\mathrm{Jac}(\phi^{-1})(y)|\,dy,
\]
for some $g\in H$. Since it is in one--to--one correspondence with $H$, the dimension of $H_*$ is also $n$.
Now, by the change of variables formula this last equals the squared $\text{L}^2$ norm of $g$ which is finite since $H\subseteq\mathrm{L}^2(\M_1,\C)$.

\noindent We now prove the formula for $K_{H*}$. Let $\varphi_1,\ldots,\varphi_n$ be an orthonormal basis of $H$. Then, $\varphi_{i,*}=\varphi_i \circ \phi^{-1}(\cdot)\sqrt{|\text{\rm Jac}(\phi^{-1})(\cdot)|}$, $1\leq i\leq n$, are elements in $H_*$ and using the change of variables formula we have:
\begin{align*}
\int_{y\in \M_2} \varphi_{i,*}(y)\overline{\varphi_{j,*}(y)}\,dy=&\int_{y\in \M_2} \varphi_i \circ \phi^{-1}(y)\,\overline{\varphi_j \circ \phi^{-1}(y)}\,|\text{\rm Jac}(\phi^{-1})(y)|\,dy\\
=&\int_{x\in\M_1}\varphi_i(x)\overline{\varphi_j(y)}\,dx=\delta_{ij},
\end{align*}
where we use the Kronecker delta notation. Hence, $\{\varphi_{i,*}\}$ form an orthonormal basis and
\begin{align*}
 K_{H_*}(a,b)=&\sum_{i=1}^n\varphi_{i,*}(a)\overline{\varphi_{i,*}(b)}\\
 =&\sum_{i=1}^n\varphi_i \circ \phi^{-1}(a)\overline{\varphi_i \circ \phi^{-1}(b)}\sqrt{|\text{\rm Jac}(\phi^{-1})(a)\,\text{\rm Jac}(\phi^{-1})(b)|}\\
 =&K_H(\phi^{-1}(a), \phi^{-1}(b))\sqrt{|\text{\rm Jac}(\varphi^{-1}(a))\, \text{\rm Jac}(\phi^{-1}(b))|}.
\end{align*}
The other formula for $K_{H_*}$ follows from this last one, using that
\[
 \mathrm{Jac}(\phi)(\phi^{-1}(a))=\mathrm{Jac}(\phi^{-1})(a)^{-1}.
\]

$\square$

\subsection{Proof of Lemma \ref{propdef_1}}\label{propLemma32}
From Proposition \ref{prop_1}, $\mathfrak{X}_*^{(r,d)}$ has reproducing kernel
\begin{equation*}
K_*^{(r,d)}(p,q) 
= 
\frac{K(\psi_d^{-1}(p), \psi_d^{-1}(q))}{\sqrt{|\text{\rm Jac}(\psi_d)(\psi_d^{-1}(p)) \text{\rm Jac}(\psi_d)(\psi_d^{-1}(q))|}}.
\end{equation*}
The Jacobian of $\psi_d$ is: 
\begin{equation}\label{eq:jacobianpsid}
\text{\rm Jac}(\psi_d)(z)
=
\left( \frac{1}{1 + ||z||^{2}} \right)^{d+1}.
\end{equation}
We thus have (denoting $p=(z,1)$ and $q=(w,1)$):
\begin{align*}
\left| K_*^{(r,d)}(p,q) \right|
=&
\frac{
\frac{r d!}{\pi^{d}}
\frac{
\left| 1 + \left\langle \psi^{-1}(p), \psi^{-1}(q) \right\rangle \right|^{L}
}{
\left( 1 + ||\psi^{-1}(p)||^{2} \right)^{\frac{L+d+1}{2}}
\left( 1 + ||\psi^{-1}(q)||^{2} \right)^{\frac{L+d+1}{2}}
}
}{
\left( \frac{1}{1 + ||\psi^{-1}(p)||^{2}} \right)^{\frac{d+1}{2}}
\left( \frac{1}{1 + ||\psi^{-1}(q)||^{2}} \right)^{\frac{d+1}{2}}
}\\
=&
\frac{r d!}{\pi^{d}} \frac{\left| 1 + \left\langle z,w \right\rangle \right|^{L}}{(1 + \|z\|^{2})^{\frac{L}{2}} (1 + \|w\|^{2})^{\frac{L}{2}}}\\
=&\frac{r d!}{\pi^{d}} \frac{\left|\left\langle p,q \right\rangle\right|^{L}}{\|p\|^{L} \|q\|^{L}},
\end{align*}
and the lemma follows.

$\square$

\subsection{Proof of Theorem \ref{Prop_proy}}\label{sec:proofofProp_proy}

Let $J$ be the quantity we want to compute. Following Proposition \ref{prop:MS} we have that
\begin{equation*}
\begin{split}
J=&  \Esp{\Xr}{ \sum_{i \neq j} \frac{1}{\sin \left( d_{\P} (x_{i}, x_{j}) \right)^{s}}  } \\&
=  \int_{\P \times \P} \frac{K(p,p)^{2} - |K(p,q)|^{2}}{\sin \left( d_{\P} (p, q) \right)^{s}} dp dq
\\
& = \frac{r^{2} d!^{2} }{\pi^{2d}}\int_{\P \times \P} \frac{1 -  \left| \left\langle p,q \right\rangle \right|^{2L} }{\left(1 - \left| \left\langle p,q \right\rangle \right|^{2}\right)^{\frac{s}{2}}} dp dq,
\end{split}
\end{equation*}
where we choose unit norm representatives $p,q$. Since the integrand only depends on the distance between $p$ and $q$ and $\P$ is a homogeneous space, we can fix $p=e_1=(1,0,\ldots,0)$ to get:
\[
 J=\frac{r^{2} d! }{\pi^{d}} \int_{\P} \frac{ 1 - \left| \left\langle e_{1},q \right\rangle \right|^{2L} }{\left(1 - \left| \left\langle e_{1},q \right\rangle \right|^{2}\right)^{\frac{s}{2}}} dq,
\]
where we have used that the volume of $\P$ is equal to ${\pi^{d}}/{d! }$. In order to compute this integral, we use the change of variables theorem with the map $\psi_d$ whose Jacobian is given in \eqref{eq:jacobianpsid}, getting:
\begin{equation*}
\begin{split}
J& =  \frac{r^{2} d! }{\pi^{d}} \int_{\C^{d}} \frac{ 1 - \left| \left\langle e_{1},\frac{(1,z)}{\sqrt{1 + ||z||^{2}}} \right\rangle \right|^{2L} }{\left(1 - \left| \left\langle e_{1},\frac{(1,z)}{\sqrt{1 + ||z||^{2}}} \right\rangle \right|^{2}\right)^{\frac{s}{2}}} \frac{1}{(1 + ||z||^{2})^{d+1}} dz
\\
& =  \frac{r^{2} d! }{\pi^{d}} \int_{\C^d} \frac{ 1 - \left( \frac{1}{1 + ||z||^{2}} \right)^{L} }{\left(1 - \frac{1}{1 + ||z||^{2}}\right)^{\frac{s}{2}}} \left(\frac{1}{1 + ||z||^{2}}\right)^{d+1} dz.
 \\
\end{split}
\end{equation*}
Integrating in polar coordinates,

\begin{equation*}
\begin{split}
J&= \frac{r^{2} d! }{\pi^{d}} \frac{2\pi^{d}}{(d-1)!} \int_{0}^{\infty} \frac{ 1 - \left( \frac{1}{1 + t^{2}} \right)^{L} }{\left(1 - \frac{1}{1 + t^{2}}\right)^{\frac{s}{2}}} \left(\frac{1}{1 + t^{2}}\right)^{d+1} t^{2d-1} dt
 \\
& =  2r^{2}d \left[ \int_{0}^{\infty} \frac{ t^{2d-1-s} }{(1 + t^{2})^{d + 1 - \frac{s}{2}}} dt - \int_{0}^{\infty} \frac{ t^{2d-1-s} }{(1 + t^{2})^{d + 1 - \frac{s}{2} + L}} dt  \right]
 \\
& =  2r^{2}d \left[ \frac{B \left( d - \frac{s}{2}, 1 \right)}{2} - \frac{B \left( d - \frac{s}{2}, L + 1 \right)}{2}  \right]
 \\
& = \frac{d}{d - \frac{s}{2}}r^{2} - r^{2}dB \left( d - \frac{s}{2}, L + 1 \right),
\end{split}
\end{equation*}
as claimed. For the assymptotics, note that for $L\to\infty$ (equiv. $r\to\infty$)
\[
 B \left( d - \frac{s}{2}, L + 1 \right)=\frac{\Gamma\left(d - \frac{s}{2}\right)\Gamma(L+1)}{\Gamma\left(d - \frac{s}{2}+ L + 1 \right)}\sim\Gamma\left(d - \frac{s}{2}\right)L^{\frac{s}{2}-d},\quad r=\binom{L+d}{d}\sim\frac{L^d}{d!},
\]
and hence
\[
 B \left( d - \frac{s}{2}, L + 1 \right)\sim\Gamma\left(d - \frac{s}{2}\right)(d!r)^{\frac{s}{2d}-1}.
\]
The assymptotic expansion claimed in the theorem follows.

$\square$
\subsection{Proof of Corollary \ref{cor:logaritmica} }\label{sec:prooflog}
Note that $\mathcal{E}_0(\omega_r)= \left. \frac{d}{ds}\right|_{s=0}\mathcal{E}_s(\omega_r)$. In particular, interchanging the order of expected value and derivative (it is an exercise to check that this change is justified), from Theorem \ref{Prop_proy} we have
\[
  \Esp{\mathfrak{X}_{*}^{(r)}}{\mathcal{E}_{0}^{\mathbb{P}}(\omega_{r}) } = \left. \frac{d}{ds} \right|_{s=0}\left(\frac{d}{d - \frac{s}{2}}r^{2} - r^{2}dB \left( d - \frac{s}{2}, L + 1 \right)\right).
\]
The proof of the corollary is now a straightforward computation of that derivative and it is left to the reader. It is helpful to recall the derivative of Euler's Beta function in terms of the digamma function $\psi_0$ for $m \in \mathbb{N}$:
\begin{align*}
 \frac{d}{dt}B(t,m)=&\frac{d}{dt}\frac{\Gamma(t)\Gamma(m)}{\Gamma(t+m)}\\=&\frac{\Gamma'(t)\Gamma(m)\Gamma(t+m)-\Gamma(t)\Gamma(m)\Gamma'(t+m)}{\Gamma(t+m)^2}\\=&\frac{\psi_0(t)\Gamma(t)\Gamma(m)-\Gamma(t)\Gamma(m)\psi_0(t+m)}{\Gamma(t+m)}\\=&B(t,m)(\psi_0(t)-\psi_0(t+m))\\=&-B(t,m)\sum_{j=0}^{m-1}\frac{1}{t+j}.
\end{align*}

\subsection{Proof of Proposition \ref{Uju1}}\label{sec:proofUju1}

We will use the following equality, valid for $y\in(-1,1)$:
\begin{equation}\label{eq:integrallibro}
\int_{0}^{2\pi} \frac{d\theta}{1 - y \cos(\theta)}=\frac{2\pi}{\sqrt{1-y^2}}.
\end{equation}
See for example \cite[3.792--1]{zwillinger2014table} from which the equality above easily follows.

We have to compute

\begin{equation}\label{eq2}
\begin{split}
& \frac{1}{(2\pi)^{r}} \int_{\theta_{1},...,\theta_{r} \in [0, 2\pi]} \Esp{\Xr}{\E{2}(y_{1}^{0},...,y_{1}^{k-1},y_{2}^{0},...,y_{r}^{k-1})} d(\theta_{1},...,\theta_{r}) 
\\
& = \frac{1}{(2\pi)^{r}} \int_{\theta_{1},...,\theta_{r} \in [0, 2\pi]} \Esp{\Xr}{ \sum_{i_{1} \neq i_{2}\text{ or } j_{1} \neq j_{2}} \frac{1}{\left|\left| y_{i_{1}}^{j_{1}} - y_{i_{2}}^{j_{2}}  \right|\right|^ {2}} } d(\theta_{1},...,\theta_{r}) 
\\
&=J_1+J_2,
\end{split}
\end{equation}
where
\begin{align*}
 J_1=& \frac{1}{(2\pi)^{r}} \int_{\theta_{1},...,\theta_{r} \in [0, 2\pi]} \Esp{\Xr}{ \sum_{i=1}^{r} \sum_{j_{1} \neq j_{2}} \frac{1}{\left|\left| y_{i}^{j_{1}} - y_{i}^{j_{2}}  \right|\right|^ {2}}}\,d(\theta_{1},...,\theta_{r}) ,
 \\
J_2=&  \frac{1}{(2\pi)^{r}} \int_{\theta_{1},...,\theta_{r} \in [0, 2\pi]}\Esp{\Xr}{ \sum_{j_{1},j_2=0}^{k-1} \sum_{i_{1} \neq i_{2}} \frac{1}{\left|\left| y_{i_{1}}^{j_{1}} - y_{i_{2}}^{j_{2}}  \right|\right|^ {2}} }\, d(\theta_{1},...,\theta_{r}).
\end{align*}
From \eqref{eq:yes} we have:
\begin{align*}
 J_1=& \frac{1}{2\pi}\sum_{i=1}^{r}\int_{\theta\in [0, 2\pi]} \Esp{\Xr}{  \sum_{j_{1} \neq j_{2}} \frac{1}{\left|\left| e^{\i\left( \theta + \frac{2\pi j_1}{k} \right)} x_i - e^{\i\left( \theta + \frac{2\pi j_2}{k} \right)} x_i   \right|\right|^ {2}}}\,d\theta.
\end{align*}
Now, the integral does not depend on $\theta$ nor in the (unit norm) vector $x_i\in\C^{n+1}$, so we actually have that
\[
 \frac{J_1}{r}=\sum_{j_{1} \neq j_{2}} \frac{1}{\left|e^{\i  \frac{2\pi j_1}{k} }  - e^{\i \frac{2\pi j_2}{k}}   \right|^ {2}},
\]
is the $2$--energy of the $k$ roots of unity. This quantity has been studied with much more detail than we need in \cite[Theorem 1.1]{BLMS:BLMS0621}. In particular, we know that it is of the form $k^3/12+o(k)$. We thus conclude:
\begin{equation}\label{eq:J1}
J_1=\frac{rk^3}{12}+o(rk).
\end{equation}
We now compute $J_2$. Interchanging the order of integration we have:


\begin{equation*}
\begin{split}
J_2
& =   \Esp{\Xr}{ \sum_{j_{1},j_2=0}^{k-1}  \sum_{i_{1} \neq i_{2}} \frac{1}{4\pi^{2}}\int_{0}^{2\pi} \int_{0}^{2\pi} \frac{d\theta_{i_{1}} d\theta_{i_{2}} }{\left|\left| e^{\i\left( \theta_{i_{1}} + \frac{2\pi j_{1}}{k} \right)} x_{i_{1}} - e^{\i\left( \theta_{i_{2}} + \frac{2\pi j_{2}}{k} \right)} x_{i_{2}}  \right|\right|^ {2}} },
\end{split}
\end{equation*}
where we can choose whatever unit norm representatives we wish of $x_{i_1}$ and $x_{i_2}$. In order to compute the inner integral, for any fixed $i_i,i_2$ we assume that our choice satisfies $\langle x_{i_1},x_{i_2}\rangle\in[0,1]$ (i.e. it is real and non--negative), which readily implies
\begin{equation}\label{eq:real}
 \sin d_{\P}(x_{i_1},x_{i_2}) 
	=
	\sqrt{1 -\left\langle x_{i_1},x_{i_2}\right\rangle ^{2}}.
\end{equation}
A simple computation using the invariance of the integral under rotations yields:
\begin{multline*}
 \frac{1}{4\pi^{2}}\int_{0}^{2\pi} \int_{0}^{2\pi} \frac{d\theta_{i_{1}} d\theta_{i_{2}} }{\left|\left| e^{\i\left( \theta_{i_{1}} + \frac{2\pi j_{1}}{k} \right)} x_{i_{1}} - e^{\i\left( \theta_{i_{2}} + \frac{2\pi j_{2}}{k} \right)} x_{i_{2}}  \right|\right|^ {2}}\\=\frac{1}{2\pi}\int_0^{2\pi}\frac{d\theta}{2-2\left\langle x_{i_1},x_{i_2}\right\rangle\cos\theta}\stackrel{\eqref{eq:integrallibro}}{=}\frac{1}{2\sqrt{1-\left\langle x_{i_1},x_{i_2}\right\rangle^2}}\stackrel{\eqref{eq:real}}{=}\frac{1}{2\sin d_{\P}(x_{i_1},x_{i_2})},
\end{multline*}
and this last value is independent of $j_1,j_2$. We thus have:
\[
 J_2= \frac{k^2}{2}\Esp{\Xr}{  \sum_{i_{1} \neq i_{2}} \frac{1}{\sin d_{\P}(x_{i_1},x_{i_2})}}
\]
This last expected value has been computed in Theorem \ref{Prop_proy}, which yields:
\begin{equation}\label{eq:J2}
J_2=
\frac{d}{2d-1}  (kr)^{2}-\frac{d\Gamma\left(d-\frac12\right)}{2(d!)^{1-\frac1{2d}}} k^2r^{1+\frac{1}{2d}}+o\left( k^2r^{1+\frac{1}{2d}}\right). 
\end{equation}
Proposition \ref{Uju1} follows from \eqref{eq2}, \eqref{eq:J1} and \eqref{eq:J2}.

$\square$

\subsection{Proof of Theorem \ref{th:main}}

Fix $d\geq1$ and let
\[
 f(A)=\frac{A^{2 - \frac{2}{2d+1}} }{12} 
-\frac{d \Gamma\left( d - \frac{1}{2} \right) A^{1 - \frac{2}{2d+1}} }{2 \left(d!\right)^{1 - \frac{1}{2d}} }
\]
be the coefficient of $n^{1+\frac{2}{2d+1}}$ in \eqref{eq3bis2}. The function $f(A)$ has a strict global minimum at
\[
 A_d=\frac{3 \Gamma\left( d - \frac{1}{2} \right) (2d-1)}{2 \left(d!\right)^{1 - \frac{1}{2d}} }.
\]
Indeed,
\[
 f(A_d)=-\frac{3^{1 - \frac{2}{2d+1}} (2d-1)^{1 - \frac{2}{2d+1}} (2d+1) \Gamma\left( d- \frac{1}{2} \right)^{2 - \frac{2}{2d+1}}}{2^{4 - \frac{2}{2d+1}} (d!)^{2 - \frac{4}{2d+1}}},
\]
gives the bound for the $\limsup$ given in Theorem \ref{th:main}. We cannot just let $k=A_dr^{\frac{1}{2d}}$ in Corollary \ref{cor:casi} since it might happen that $k\not\in\mathbb{Z}$, but we will easily go over this problem. Let $L\geq1$ be any positive integer, let $r=\binom{d+L}{d}$ and let $A$ be the unique number in the interval
\[
 \big[A_d,A_d+r^{-\frac{1}{2d}}\big)
\]
such that $k=Ar^{\frac{1}{2d}}\in\mathbb{Z}$. Finally, let $n=n_L=rk$, which depends uniquely on $d$ and $L$, and which satisfies $n_L\to\infty$ as $L\to\infty$. For any $\epsilon>0$ we then have:
\begin{align*}
 \limsup_{L\to\infty}\frac{ V_{2}(\mathbb{S}^{2d+1})n_L^{2}-\min_{\omega_{n_L}}\left(\E{2}(\omega_{n_L})\right) }{n_L^{1+\frac{2}{2d+1}}}
\geq -f(A)\geq -f(A_d)-\epsilon,
\end{align*}
the first inequality from Corollary \ref{cor:casi} and the second inequality due to $r\to\infty$ as $L\to\infty$, which implies for some constant $C>0$:
\[
|f(A)-f(A_d)| \leq Cr^{-\frac{1}{2d}}\to0,\quad L\to\infty.
\]
We have thus proved
\[
  \limsup_{L\to\infty}\frac{ V_{2}(\mathbb{S}^{2d+1})n_L^{2}-\min_{\omega_{n_L}}\left(\E{2}(\omega_{n_L})\right) }{n_L^{1+\frac{2}{2d+1}}}
\geq -f(A_d),
\]
which finishes the proof of our Theorem \ref{th:main}.


\subsection{Proof of Theorem \ref{th:main2}}

From \cite{Juan}, the Green function of $\P$ is given $G(x,y)=\phi(r)$ where $r=d_{\P}(x,y)$ and
\[
 \phi'(r)=-\frac{1}{Vol(\P)}\frac{\int_r^{\pi/2}\sin^{2d-1}t\cos t\,dt}{\sin^{2d-1}r\cos r}=-\frac{1}{2dVol(\P)}\frac{1-\sin^{2d}r}{\sin^{2d-1}r\cos r}.
\]
Integrating the formula above (see for example \cite[2.517--1]{zwillinger2014table} we have:
\begin{equation*}
\begin{split}
\phi(r)
= 
& \frac{1}{2dVol(\P)}
\left[
\frac{1}{2} \displaystyle\sum_{k=1}^{d-1} \frac{1}{(d-k) \left( \sin r \right)^{2d-2k}} 
- 
\log \left( \sin r \right)
\right] +C.
\end{split}
\end{equation*}
In order to compute the constant we need to impose that the average of $G(x,\cdot)$ equals $0$ for all (i.e. for some) $x\in\P$. Let $x=(1,0)$ and change variables using $\psi_d$ from Lemma \ref{propdef_1} whose Jacobian is given in \eqref{eq:jacobianpsid} to compute:
\[
 C=- \frac{1}{2dVol(\P)^2}
\left[
\frac{1}{2} \displaystyle\sum_{k=1}^{d-1} \int_{z\in\C^d}\frac{(1+\|z\|^2)^{-k-1}}{(d-k) \|z\|^{2d-2k}} 
\,dz
- 
\frac{1}{2}\int_{z\in\C^d}\frac{\log \left( \frac{\|z\|^2}{1+\|z\|^2} \right)}{(1+\|z\|^2)^{d+1}}\,dz
\right].
\]
Integrating in polar coordinates,
\begin{align*}
 C=&
\frac{1}{2Vol(\P)}\left(\int_{0}^\infty\frac{t^{2d-1}\log \left( \frac{t^2}{1+t^2} \right)}{(1+t^2)^{d+1}}\,dt-  \displaystyle\sum_{k=1}^{d-1} \int_{0}^\infty\frac{t^{2k-1}}{(d-k)(1+t^2)^{k+1}} \,dt\right)
\\=&
-\frac{d!}{4\pi^d}\left(\frac{1}{d^{2}}
+
\displaystyle\sum_{k=1}^{d-1} \frac{1}{k(d-k)}\right)\\=&
-\frac{(d-1)!}{4\pi^d}\left(\frac{1}{d}+2\sum_{k=1}^{d-1}\frac1k\right).
\end{align*}
(for the computation of the integrals, use the change of variables $s=t^2/(1+t^2)$ and \cite[4.272--15]{zwillinger2014table}, for example).

We thus conclude for $r=d_{\P}(x,y)$ :
\begin{equation*}
\begin{split}
G(x,y) 
= 
& \frac{(d-1)!}{2\pi^{d}}
\left[
\left(
\frac{1}{2} \displaystyle\sum_{k=1}^{d-1} \frac{1}{(d-k) \left( \sin r \right)^{2d-2k}} 
\right)
- 
\log \left( \sin r\right)
\right] \\
& 
-\frac{(d-1)!}{4\pi^d}\left(\frac{1}{d}+2\sum_{k=1}^{d-1}\frac1k\right).
\end{split}
\end{equation*}

\noindent Following the definitions of Theorem \ref{Prop_proy} and Corollary \ref{cor:logaritmica}, the expected value of Green energy may be expressed as

\begin{equation*}
\begin{split}
& \Esp{\mathfrak{X}_{*}^{(r)}}{\E{G}^{\mathbb{P}}(\omega_{r})}
 =
\overbrace{\frac{(d-1)!}{2\pi^{d}}
\left[
\left(
\frac{1}{2} \displaystyle\sum_{k=1}^{d-1} \frac{1}{d-k} 
\Esp{\mathfrak{X}_{*}^{(r)}}{\mathcal{E}_{2d-2k}^{\mathbb{P}}(\omega) }
\right)
+
\Esp{\mathfrak{X}_{*}^{(r)}}{\mathcal{E}_{0}^{\mathbb{P}}(\omega) }
\right]}^{A}  \\
&-\frac{r(r-1)(d-1)!}{4\pi^d}\left(\frac{1}{d}+2\sum_{k=1}^{d-1}\frac1k\right). \\
\end{split}
\end{equation*}
Each of the expected values in the last expression has been computed in Theorem \ref{Prop_proy} and Corollary \ref{cor:logaritmica}, producing:

\begin{equation*}
\begin{split}
&
A
= 
\frac{(d-1)!}{4\pi^{d}}
\left[
\left(
\displaystyle\sum_{k=1}^{d-1} \frac{1}{d-k} 
\left(
\frac{d}{k}r^{2}-r^{2-\frac{k}{d}}\frac{d \Gamma\left(k\right)}{(d!)^{\frac{k}{d}}}
\right)
\right)
+
\frac{r^2}{d} + \frac{r\log r}{d}
\right]
 +o\left(r^{2-\frac{1}{d}}\right)\\
& =
r^{2} \frac{d!}{4\pi^{d}}
\left(
\sum_{k=1}^{d-1} \frac{1}{k(d-k)} + \frac{1}{d^2}
\right)
-
\frac{(d!)^{1- \frac{1}{d}}}{4\pi^{d}(d-1)}
  r^{2 - \frac{1}{d}}
+o\left(r^{2-\frac{1}{d}}\right)
\\
& =
r^{2} \frac{(d-1)!}{4\pi^{d}}
\left(
\frac{1}{d}+2\sum_{k=1}^{d-1} \frac{1}{k} 
\right)
-
\frac{(d!)^{1- \frac{1}{d}}}{4\pi^{d}(d-1)}
  r^{2 - \frac{1}{d}}
+o\left(r^{2-\frac{1}{d}}\right)
\end{split}
\end{equation*}

%
%

\noindent We thus have:

\begin{equation*}
\begin{split}
\Esp{\mathfrak{X}_{*}^{(r)}}{\E{G}^{\mathbb{P}}(\omega_{r})}
&  =
-
\frac{(d!)^{1- \frac{1}{d}}}{4\pi^{d}(d-1)}
  r^{2 - \frac{1}{d}}
+o\left(r^{2-\frac{1}{d}}\right).
\end{split}
\end{equation*}

Since this last equation holds for an infinite sequence of numbers (those of the form $r=\binom{d+L}{d}$, Theorem \ref{th:main2} follows.



\end{document}